\documentclass[a4paper, 12pt]{article}
\usepackage{amssymb,latexsym}
\usepackage[cp1251]{inputenc}
\usepackage[english]{babel}
\usepackage{fullpage,amsfonts,amsthm}
\usepackage[intlimits]{amsmath}
\usepackage{indentfirst}
\newtheorem{theorem}{Theorem}

\title {
\vspace{15mm}
\textbf{On an equation involving fractional powers with one prime and one almost prime variables}}
\author{
\vspace{8mm}
Zh. H. Petrov \hspace{15mm} D. I. Tolev}
\date{}

\begin{document}

\maketitle
\begin{abstract}
In this paper we consider the equation $[p^{c}] + [m^{c}] = N$, where $N$ is a sufficiently large integer,
and prove that if $1 < c < \frac{29}{28}$, then it has a solution in a prime $p$
and an almost prime $m$ with at most $\left[ \frac{52}{29 - 28 c}\right] + 1 $
prime factors.
\end{abstract}

\section{Introduction and statement of the result}

A Piatetski--Shapiro sequence is a sequence of the form
\begin{equation} \label{1}
\{ [n^{c}]\}_{n \in \mathbb N } \qquad (c > 1, c \not\in \mathbb{N}),
\end{equation}
where $[t]$ denotes the integer part of $t$. In 1953 Piatetski--Shapiro~\cite{12} showed that if
$1 < c < \frac{12}{11}$, then the sequence \eqref{1} contains infinitely many prime numbers.
Since then, the upper bound for $c$ has been improved many times and
the strongest result is due to Rivat and Wu~\cite{Rivat_Wu}. They proved that
the sequence \eqref{1} contains infinitely many primes provided that
$1 < c < \frac{243}{205}$.

\bigskip

For any natural number $r$, let $\mathcal P_{r}$ denote the set of $r$-almost primes,
i.e. the set of natural numbers having at most $r$ prime factors counted with multiplicity.
There are many papers devoted to the study of
problems involving Piatetsi--Shapiro primes and almost primes.
In 2011 Cai and Wang~\cite{Cai_Wang}, improving an earlier result of Peneva~\cite{Peneva}, 
showed that if $1 < c < \frac{30}{29}$, then there exist infinitely many primes $p$ of the form 
$[n^c]$ such that $p + 2 \in \mathcal P_5$.
Later, in 2014, Baker, Banks, Guo, and Yeager~\cite{10} showed that if
$1 < c < \frac{77}{76}$, then the sequence
\begin{equation*}
\{ [n^{c}]\}_{n \in \mathcal P_{8}} 
\end{equation*}
contains infinitely many prime numbers.

\bigskip

Consider the equation
\begin{equation}\label{2}
[m_{1}^{c}] + [m_{2}^{c}] = N .
\end{equation}
In 1973, Deshouillers~\cite{1} proved that if $1 < c < \frac{4}{3}$, 
then for every sufficiently large integer $N$ the equation \eqref{2} has a solution with $m_1$ and $m_2$ integers.
This result was improved by Gritsenko~\cite{6}, and later by Konyagin~\cite{7}.
In particular, the latter author showed that \eqref{2} has a solution in integers $m_{1}$, $m_{2}$
for $1 < c < \frac{3}{2}$ and $N$ sufficiently large.

\bigskip

Kumchev~\cite{2} proved that if $1 < c < \frac{16}{15}$, then every sufficiently large integer $N$
can be represented in the form \eqref{2}, where $m_{1}$ is a prime and $m_{2}$ is an integer.
On the other hand, the celebrated theorem of Chen~\cite{Chen} states that every sufficiently large even integer
can be represented as a sum of a prime and an almost prime from $\mathcal P_2$.
Having in mind this profound result, one can conjecture that there exists a constant $c_0 > 1$ such that
if $1 < c < c_0$, then the equation \eqref{2} has a solution 
with $m_1$ a prime and $m_2 \in \mathcal P_2$ provided that $N$ is sufficiently large.
In the present paper, we establish a result of this type and
prove the following
\begin{theorem} 
Suppose that $1 < c < \frac{29}{28}$. Then every sufficiently large integer $N$ can be represented as
\begin{equation} \label{1001}
  [p^c] + [m^c] = N ,
\end{equation}
where $p$ is a prime and $m$ is an almost prime with at most
$\left [ \frac{52}{29-28c} \right ] + 1$ prime factors.
\end{theorem}

\bigskip

We note that the integer $\left [ \frac{52}{29-28c} \right ] + 1$ is equal to $53$ if $c$ is close to $1$
and it is large if $c$ is close to $\frac{29}{28}$.

\bigskip

Our first step in the proof is to apply the linear sieve. After doing so, we could try to establish a relatively
strong estimate for the exponential sum defined in \eqref{24} which is a rather difficult task since the function in the exponent depends on $[p^c]$. Instead, we represent this sum as a linear combination of similar sums 
(see \eqref{39}) with a smooth function of $p$ in the exponent. Then we use standard techniques to estimate these sums. We would like to mention that the sums in \eqref{39} are also studied by Kumchev~\cite{2}. However, we cannot use his work because we require stronger bounds for them.

\section{Notation}

We fix the following notation: $\{ t \}$ is the fractional part of $t$,  
the function $\rho(t)$ is defined by 
$\rho(t) = \frac{1}{2} - \{ t \}$ and $e(t) = e^{2\pi i t}$.
We use Vinogradov's notation $A \ll B$, which is equivalent to $A = O(B)$. If we have
simultaneously $A \ll B$ and $B \ll A$, then we shall write $A \asymp B$.

\bigskip

For us $p$ will be reserved for prime numbers.
By $\varepsilon$ we denote an arbitrarily small positive number,
which is not necessarily the same in the different formulae. 
As usual, $\sum_{n \le x}$ means $\sum_{1 \le n \le x}$ and $\mu(n)$, $\Lambda(n)$ and $\tau(n)$ are
the Mobius function, von Mangolds' function and the number of positive divisors of $n$, respectively.

\section{Proof of the theorem}

\subsection{Beginning of the proof}

Let $N$ be a sufficiently large integer and let
\begin{equation}\label{3}
1 < c < \frac{29}{28},\hspace{15mm} \gamma = \frac{1}{c},\hspace{15mm} P = 10^{-9} N^{\gamma}.
\end{equation}
Suppose that $\alpha > 0$ is a constant, which will be specified later, and let
\begin{equation}\label{4}
z=N^{\alpha }, \hspace{15mm}   B_{z} = \prod_{p < z}{p}.
\end{equation}
We consider the sum
\begin{equation}
\Gamma = \sum_{\substack{P < p \leq 2P , \; m \in \mathbb{N} \\ [p^c] + [m^c] = N \\ (m, B_{z}) = 1}}{\log p}. \label{5}
\end{equation}

If $\Gamma > 0$, then there is a prime number $p$ and a natural number $m$ satisfying
the conditions imposed in the domain of summation of $\Gamma $. 
From the condition $(m, B_{z})=1$ it follows that any prime factor of $m$
is greater or equal to $z$. 
Suppose that $m$ has $l$ prime factors, counted with the multiplicity. 
Then we have
\[
N^{\gamma } \geq m  \geq z^{l} = N^{\alpha l} 
\]
and hence $l \leq \frac{\gamma }{\alpha }$.
This implies that if $\Gamma > 0$ then \eqref{1001} has a solution 
with $p$ a prime and $m$ an almost prime with at most
$\left [ \frac{\gamma }{\alpha } \right ]$ prime factors.

\bigskip

We denote
\begin{equation}\label{6}
D=N^{\delta },
\end{equation}
where $\delta > 0$ is a constant which will be specified later.
Let
$\lambda (d)$ be the lower bound Rosser weights of level $D$, 
(see \cite[Chapter~$4$]{8}).
Then we have
\begin{equation} \label{7}
\sum_{d|k}{\mu (d)} \geq \sum_{d|k}{\lambda (d)}
 \qquad \text{for every} \qquad  k \in \mathbb{N} .
\end{equation}
Furthermore, we know that
\begin{equation} \label{8}
 |\lambda (d)| \leq 1 \;\; \text{ for all } \; d \, ; \qquad
 \lambda (d) = 0 
   \;\;  \text{ if } \;\; d > D \;\; \text{ or } \;\; \mu(d) = 0 .
\end{equation}
Finally, we have
\begin{equation} \label{9}
\sum_{d \mid B_{z}}{\frac{\lambda (d)}{d}} \geq  \prod_{p < z} \left( 1 - \frac{1}{p} \right)
\left ( f(s) + O ((\log D)^{-\frac{1}{3}})\right ),
\end{equation}
where
\begin{equation} \label{10010}
  s = \frac{\log D}{\log z} = \frac{\delta }{\alpha }
\end{equation}
and where $f(s)$ is the lower function of the linear sieve, for which we know that
\begin{equation} \label{10020}
 f(s) = 0 \;\;\; \text{for} \;\;\; 0 < s < 2 ;
 \qquad
 f(s) = 2e^{G} s^{-1}\log (s-1) \;\;\; \text{for} \;\;\;  2<s<3 .
\end{equation}
(Here $G$ is the Euler constant).

\bigskip

From \eqref{5} and \eqref{7} we find
\begin{equation*}
\Gamma = \sum_{\substack {P < p \leq 2P , \; m \in \mathbb N \\ [p^c] + [m^c] = N}}
{(\log p)}\sum_{d \mid (m, B_{z})}{\mu (d)}
\geq \sum_{\substack {P < p \leq 2P , \; m \in \mathbb N  \\ [p^c] + [m^c] = N}}{(\log p)}\sum_{d|(m, B_{z})}
{\lambda (d)}.
\end{equation*}
We change the order of summation to obtain
\begin{equation}\label{10}
\Gamma \geq \sum_{d \mid B_{z}}{\lambda (d) \, G_{d}} \, , \qquad \text{where} \hspace{10mm }
G_{d} = \sum_{\substack {P < p \leq 2P , \; 
m \in \mathbb N \\ [p^c] + [m^c] = N \\ m \equiv 0(\bmod {d}) }}{\log p}.
\end{equation}
Now, we write the sum $G_{d}$ in the form
\begin{equation} \label{11}
G_{d} = \sum_{P < p \leq 2P}{(\log p) \, G_{d, p}'} \, , \qquad \text{where} \hspace{10mm }
G_{d, p}' = \sum_{\substack {m \in \mathbb N \\ m \equiv 0(\bmod d) \\ [p^c] + [m^c] = N}}{1}.
\end{equation}
We use the obvious identity
\[
  \sum_{a \leq m < b}{1} = [-a] - [-b] = b - a - \rho (-b) + \rho (-a)
\]
to establish
\begin{align}
G_{d, p}'
  & =
  \sum_{\substack {m \in \mathbb N \\ m \equiv 0(\bmod {d}) \\ N-[p^c]\leq m^c < N+1-[p^c]}}{1}
     = \sum_{ \frac{1}{d}(N - [p^c])^{\gamma } \le m < \frac{1}{d}(N + 1 - [p^c])^{\gamma }} {1}
   \notag \\
   & \notag \\
  & =
    \frac{(N + 1 - [p^c])^{\gamma } - (N - [p^c])^{\gamma }}{d}
       + \rho \left ( -\frac{1}{d}(N - [p^c])^{\gamma}\right )
       - \rho \left ( -\frac{1}{d}(N + 1 - [p^c])^{\gamma }\right ).
       \notag
\end{align}
We combine the above with \eqref{11} to obtain
\begin{equation}\label{12}
G_{d} = \frac{1}{d}A(N) + \sum_{P<p\leq 2P}{(\log p)\left ( \rho \left ( -\frac{1}{d}(N - [p^c])^{\gamma}\right ) -
\rho \left ( -\frac{1}{d}(N + 1 - [p^c])^{\gamma}\right ) \right )},
\end{equation}
where
\begin{equation*}
A(N) = \sum_{P<p\leq 2P}(\log p){\left ((N-[p^c]+1)^{\gamma } - (N-[p^c])^{\gamma }\right )}.
\end{equation*}

\bigskip

From  
\begin{equation*}
(N-[p^c]+1)^{\gamma } = (N - [p^c])^{\gamma } + \gamma (N - [p^c])^{\gamma - 1} + O(N^{\gamma - 2}) ,
\end{equation*}
we deduce that 
\begin{equation*}
  A(N) = \gamma \sum_{P < p \leq 2P}{(\log p) \left ( (N-[p^c])^{\gamma -1 }+O(N^{\gamma -2})\right )},
\end{equation*}
and by Chebyshev's prime number theorem and the definition of $P$ in \eqref{3}, we get
\begin{equation}\label{14}
A(N) \asymp N^{2\gamma -1}.
\end{equation}

\bigskip

From \eqref{10} and \eqref{12} we have
\begin{equation}\label{15}
\Gamma \geq \Gamma_{0} + \Sigma_{0} - \Sigma_{1},
\end{equation}
where
\begin{align}
  & 
  \Gamma_{0} = A(N)\sum_{d \mid B_{z}}{\frac{\lambda (d)}{d}}, 
  \label{16}  \\
  & 
  \Sigma_{j} = \sum_{d \mid B_{z}}{\lambda (d)\sum_{P<p\leq 2P}{(\log p) \, \rho
  \left ( -\frac{1}{d}(N + j - [p^c])^{\gamma}\right )}}, \hspace{5mm} j = 0, 1. 
  \label{17}
\end{align}

\bigskip

Consider $\Gamma_0$.
We use \eqref{4} and the Mertens formula to find
\begin{equation} \label{10025}
 \prod_{p < z}{\left ( 1 - \frac{1}{p}\right )} \asymp \frac{1}{\log z} \asymp \frac{1}{\log N}.
\end{equation}
Assume that 
\begin{equation} \label{10030}
   2 < \frac{\delta }{\alpha} < 3 .
\end{equation}
Then, having in mind \eqref{10010} and \eqref{10020} we find that $f(s) > \kappa$ for some constant
$\kappa > 0$ depending on $\delta$ and $\alpha$ only.
Therefore, using  \eqref{9} and \eqref{10025}
we get
\[
 \sum_{d \mid B_{z}}{\frac{\lambda (d)}{d}}
 \gg \frac{1}{\log N} .
\]
Thus, by \eqref{14} and \eqref{16} we obtain
\begin{equation}\label{18}
\Gamma_{0} \gg \frac{N^{2\gamma -1}}{\log N}.
\end{equation}

\bigskip

We aim to establish the following bound
for the sums $\Sigma_{j}$ defined in \eqref{17}:
\begin{equation} \label{10040}
  \Sigma_{j}  \ll \frac{N^{2\gamma - 1}}{(\log N)^2} , \qquad j = 0, 1 .
\end{equation}
This, together with \eqref{15} and \eqref{18} would imply
\begin{equation*}
\Gamma \gg \frac{N^{2\gamma -1}}{\log N},
\end{equation*}
hence $\Gamma > 0$ for sufficiently large $N$.
Then, as we already explained, the equation \eqref{1001} would have a solution in a prime $p$
and an almost prime $m$ with no more than $\left[ \frac{\gamma}{\alpha} \right]$
prime factors.

\bigskip

The remaining part of the paper will be devoted to the proof of the estimates \eqref{10040}
under the assumptions
\begin{equation} \label{10050}
  \frac{28}{29} < \gamma < 1 , \qquad
  \delta < \frac{29 \gamma - 28}{26} .
\end{equation}
Then, it would remain to choose 
\[
  \alpha = \frac{29 \gamma - 28}{52} - \varepsilon_0
\]
for some small $\varepsilon_0 > 0$ and to take
\[
   \delta \in \left( 2 \alpha , \frac{29 \gamma - 28}{26} \right) .
\]

In this case, when $\varepsilon_0$ is small enough the condition \eqref{10030} holds. 
With the choice \eqref{3} of $c$, it is easy to see that
$\left[ \frac{\gamma}{\alpha} \right] \leq \left[ \frac{52}{29 - 28 c} \right] + 1$, which proves the theorem.

\subsection{The estimation of the sums $\Sigma_{1}$ and $\Sigma_2 $ --- beginning}

Consider the sum $\Sigma _{j}$ defined in \eqref{17}. We apply Vaaler's theorem~\cite{11},
which states that
for each $H \geq 2$ there are numbers $c_{h}$ ($0<|h|\leq H$), $d_{h} (|h| \leq H)$, such that
\begin{equation}\label{20}
\rho (t) = \sum_{0<|h|\leq H} {c_{h} \, e(ht)} + \Delta_{H}({t}),
\end{equation}
where
\begin{equation}\label{94}
|\Delta_{H}(t)| \leq \sum_{|h|\leq H} {d_{h} \, e(ht)}
\end{equation}
and
\begin{equation}\label{19}
|c_{h}| \ll \frac{1}{|h|}, \hspace{5mm} |d_{h}| \ll \frac{1}{H}.
\end{equation}
From \eqref{17} and \eqref{20} it follows that
\begin{equation}\label{21}
\Sigma_{j} = \Sigma_{j}' + \Sigma_{j}'',
\end{equation}
where
\begin{align}
\Sigma_{j}'
    & =
        \sum_{d \leq D}{\lambda (d)}\sum_{P<p\leq 2P}{(\log p)}
        \sum_{0<|h|\leq H}{c_{h} \, e\left ( -\frac{h}{d}(N + j -[p^{c}])^{\gamma }\right )}, 
        \label{22} \\
     & \notag \\   
\Sigma_{j}''
    & =
        \sum_{d \leq D}{\lambda (d)}
        \sum_{P < p \leq 2P}{(\log p) \, \Delta_{H} \left ( -\frac{(N + j - [p^{c}])^{\gamma }}{d}\right )}. \label{23}
\end{align}

\bigskip

Let
\begin{equation}\label{24}
W(v) = \sum_{P < p \leq 2P}{(\log p) \, e(v(N + j - [p^{c}])^{\gamma })} .
\end{equation}
We start with the sum $\Sigma_{j}'$. Changing the order of summation together with 
\eqref{8}, \eqref{19} and \eqref{24} implies that
\begin{equation}\label{26}
\Sigma_{j}' = \sum_{d \leq D}{\lambda (d)}\sum_{0 < |h| \leq H}
   {c_{h} \, W \left( - \frac{h}{d} \right) } \ll \sum_{d \leq D}\sum_{1 \leq h \leq H}{\frac{1}{h}
   \left | W \left( \frac{h}{d} \right) \right |},
\end{equation}
For the sum $\Sigma_{j}''$ we use 
\eqref{8}, \eqref{94}, \eqref{19}  and \eqref{24} to get
\begin{align*}
\Sigma_{j}''
    & \ll
        \sum_{d \leq D} \sum_{P < p \leq 2P}{(\log p)}\sum_{|h| \leq H}{d_{h} \,
        e\left ( -\frac{h}{d} (N + j - [p^{c}])^{\gamma }\right )} 
        \\
      & \\  
    & =
        \sum_{d \leq D} \sum_{|h| \leq H}{d_{h} \, W \left( - \frac{h}{d} \right)}
        \ll
            \sum_{d \le D} \sum_{|h| \le H} \frac{1}{H} \, \left| W \left( \frac{h}{d}\right) \right|.
\end{align*}
From \eqref{3}, \eqref{24} and Chebyshev's prime number theorem we find that $W(0) \asymp N^{\gamma }$
and hence
\begin{equation}\label{27}
\Sigma_{j}'' \ll \sum_{d \leq D}{\frac{N^{\gamma }}{H}}
+ \sum_{d \leq D}\sum_{1 \leq h \leq H}{\frac{1}{H}\left | W \left ( \frac{h}{d}\right )\right |}.
\end{equation}
We let
\begin{equation}\label{28}
  H = dN^{1 - \gamma }(\log N)^{3} .
\end{equation}
Now, using \eqref{6}, \eqref{21} and \eqref{26} -- \eqref{28} we obtain
\begin{equation}\label{83}
\Sigma_{j} \ll \frac{N^{2\gamma - 1}}{(\log N)^{2}} + \sum_{d \leq D}\sum_{h \leq H}\frac{1}{h}
   \left | W \left(  \frac{h}{d} \right) \right |, \hspace{5mm} j=0, 1.
\end{equation}

\subsection{Consideration of the sum $W(v)$}

As we mentioned earlier, it is hard to estimate directly the exponential sum $W(v)$, defined by \eqref{24}. 
Instead, we can write it as a linear combination of similar sums which are easier to be dealt with.

\bigskip

Let $ Z \ge 2$ be an integer, which we shall specify later.
We apply the well-known Vinogradov's ``little cups'' lemma
(see \cite[Chapter 1, Lemma A]{4}) with parameters
\[
  \alpha = - \frac{1}{4 Z} , \qquad \beta =  \frac{1}{4 Z} , \qquad
  \Delta =  \frac{1}{2 Z} , \qquad r = [\log N]
\]
and construct a function $g(t)$ which is periodic with period $1$ and has the following properties:
\begin{equation} \label{17.02.eq5}
  g(0) = 1 ; \quad 0 < g(t) < 1 \;\; \text{for} \;\; 0 < |t| < \frac{1}{2 Z} ;
  \quad g(t) = 0 \;\; \text{for} \;\; \frac{1}{2 Z} \le |t| \le \frac{1}{2} .
\end{equation}

\bigskip

Furthermore, the Fourier series of $g(t)$ is given by
\[
  g(t) = \frac{1}{2 Z} + \sum_{\substack{n \in \mathbb Z \\ n \not= 0}} \beta_n \, e (nt) ,
   \quad \text{with} \quad
   |\beta_n| \le \min
   \left ( \frac{1}{2 Z}, \frac{1}{|n|}\left ( \frac{2 \, Z \, [\log N] }{\pi |n|} \right)^{[\log N]}  \right) .
\]
From the above estimate of $|\beta_n|$ one easily obtains
\[
  \sum_{|n| > Z (\log N)^4 } |\beta_n| \ll N^{- \log \log N}
\]
with an absolute constant in the $\ll$-symbol.
Hence we have
\begin{equation} \label{17.02.eq10}
  g(t) = \sum_{ |n| \le Z (\log N)^4 } \beta_n \, e (nt)
   + O \left( N^{- \log \log N}  \right) ,
\end{equation}
where the implied constant is absolute and
\begin{equation} \label{08.03.2016.eq50}
  |\beta_n| \le \frac{1}{2 Z}  .
\end{equation}

\bigskip

Finally, one can easily see that the function $g(t)$, constructed in the proof of 
\cite[Chapter~1, Lemma~A]{4}, is even and also satisfies
\begin{equation} \label{17.02.eq20}
  g(t) + g \left( t - \frac{1}{2 Z} \right) = 1
  \qquad \text{for} \qquad 0 \le t \le \frac{1}{2 Z} .
\end{equation}

\bigskip

Let
\begin{equation} \label{17.02.eq30}
  g_z(t) = g \left( t - \frac{z}{2Z} \right) \qquad \text{for} \qquad
  z = 0, 1, 2, \dots, 2Z-1 .
\end{equation}
Obviously, each $g_z(t)$ is a periodic function with period $1$.
From \eqref{17.02.eq5} we find that
\begin{align}
  &
   0 < g_z(t) \le 1 \qquad \text{if} \qquad
   \left| t - \frac{z}{2Z}  \right| < \frac{1}{2Z}
    ;
   \label{17.02.eq40} \\
   \notag \\
  &
   g_z (t) = 0 \qquad \qquad \text{if} \qquad
    \frac{1}{2Z} \le \left| t - \frac{z}{2Z}  \right| \le \frac{1}{2} .
     \label{17.02.eq50}
\end{align}
From \eqref{17.02.eq30} it follows that if $\beta_n^{(z)}$ is the $n$-th Fourier coefficient of the function
$g_z(t)$, then $ \beta_n^{(z)} = \beta_n \, e \left( - \frac{zn}{2M} \right)$ and hence
$ | \beta_n^{(z)} | = |\beta_n| $.
From this observation and \eqref{17.02.eq10}, 
as well as the estimate for $|\beta_n|$ given above,
we find that for $z =0, 1, \dots , 2Z-1$ we have
\begin{equation} \label{17.02.eq60}
  g_z(t) = \sum_{  |n| \le Z (\log N)^4 } \beta_n^{(z)} \, e (nt)
   + O \left( N^{- \log \log N}  \right) ,
\end{equation}
where the constant in the $O$-symbol is absolute and
\begin{equation} \label{17.02.eq70}
 |\beta_n^{(z)}| \le \frac{1}{2Z} .
\end{equation}

\bigskip

Finally, from \eqref{17.02.eq20}, \eqref{17.02.eq30} and \eqref{17.02.eq50} we find
\begin{equation} \label{17.02.eq80}
   \sum_{z=0}^{2Z-1} g_z (t) = 1 \qquad \text{for all} \quad t \in \mathbb R .
\end{equation}
(We leave the easy verification to the reader).

\bigskip

Now, we consider the sum $W(v)$ which was defined in \eqref{24}.
From \eqref{17.02.eq80} it follows that
\begin{equation} \label{17.02.eq90}
  W(v)  =  \sum_{P < p \leq 2P} (\log p) \, e(v(N + j - [p^{c}])^{\gamma }) \,
  \sum_{z=0}^{2Z-1} g_z(p^{c})
   =  \sum_{z=0}^{2Z-1} W_z(v) ,
\end{equation}
where
\begin{equation} \label{17.02.eq100}
  W_z(v)  =  \sum_{P < p \leq 2P} (\log p) \,  g_z(p^{c}) \, e(v(N + j - [p^{c}])^{\gamma }) .
\end{equation}
It is clear that
\[
  W_0(v) \ll \sum_{P < p \leq 2P} (\log p) \,  g (p^{c}) \ll \log N \sum_{P < k \leq 2P}   g (k^{c}) .
\]
We apply \eqref{3}, \eqref{17.02.eq10} and \eqref{08.03.2016.eq50} to get
\begin{align}
 W_0(v)
   & \ll
   (\log N) \,
  \frac{P}{Z}  + (\log N) \left| \sum_{P < k \leq 2P} \; \sum_{0 < |n| \le Z (\log N)^4 } \beta_n \; e (n k^c) \right|
  + 1 
  \notag \\
  & \notag \\
   & \ll
  (\log N) \frac{N^{\gamma}}{Z} + \frac{\log N}{Z} \sum_{ n \le Z (\log N)^4 } \left| \mathcal H_n  \right| 
  + 1 ,
  \label{17.02.eq110}
\end{align}
where
\[
 \mathcal H_n =  \sum_{P < k \leq 2P}  e (n k^c) .
\]

\bigskip

If $ \theta (x) = n x^c$, then $ \theta'' (x)  = c(c-1)nx^{c-2} \asymp n P^{c-2} $
uniformly for $x \in [P, 2P]$.
Hence, we can apply Van der Corput's theorem (see \cite{4}, Chapter~1, Theorem~5) to get
\begin{equation} \label{17.02.eq120}
 \mathcal H_n \ll P \left( n P^{c-2} \right)^{\frac{1}{2}} + \left( n P^{c-2}  \right)^{- \frac{1}{2}}
 \ll P^{\frac{c}{2}} \, n^{\frac{1}{2}} .
\end{equation}

\bigskip

Henceforth we assume that
\begin{gather}\label{86}
  Z \ll N^{\frac{2\gamma -1}{3}} \, (\log N)^{-4} .
\end{gather}
Then using \eqref{3} and \eqref{17.02.eq110} -- \eqref{86}
we obtain
\begin{equation} \label{17.02.eq130}
 W_0(v) \ll (\log N) \left( \frac{N^{\gamma}}{Z} + N^{\frac{1}{2}} \, Z^{\frac{1}{2}} \, \log^6 N \right)
 \ll (\log N) \, \frac{N^{\gamma}}{Z} .
\end{equation}

Now, we restrict our attention to the sums
$W_z(v)$ for $1 \le z \le 2Z-1$.
Using \eqref{17.02.eq50} we see that $g_z(p^c)$ vanishes unless
$\{ p^c \} \in \left[ \frac{z-1}{2Z}, \frac{z+1}{2Z}\right]$.
Hence, the only summands in the sum \eqref{17.02.eq100} are those for which
$\{ p^c \} = \frac{z}{2Z} + O \left( \frac{1}{Z} \right)$.
In this case we have
\[
v(N + j - [p^{c}])^{\gamma } = v\left ( N + j - p^{c} + \frac{z}{2Z} \right )^{\gamma }
+ O \left ( \frac{vN^{\gamma -1}}{Z}\right ) \label{34}
\]
and respectively
\[
 e \left( v(N + j - [p^{c}])^{\gamma } \right)
 =
  e \left( v\left ( N + j - p^{c} + \frac{z}{2Z} \right )^{\gamma }  \right)
  + O \left ( \frac{vN^{\gamma -1}}{Z}\right ) .
\]
Then using \eqref{17.02.eq100} we find
\begin{equation} \label{08.03.2016.eq105}
  W_z(v) = V_z(v) + O \left( \frac{v N^{\gamma - 1 }}{Z}  \sum_{P < p \le 2 P } (\log p) \, g_m (p^c) \right) ,
\end{equation}
where
\begin{gather}\label{35}
V_z(v) = \sum_{P < p \leq 2P} \, (\log p) \, g_{z} (p^{c}) \,
 e \left ( v\left ( N + j - p^{c} + \frac{z}{2Z} \right )^{\gamma }\right ).
\end{gather}

\bigskip

Therefore, from \eqref{17.02.eq90}, \eqref{17.02.eq130} and \eqref{08.03.2016.eq105}
we obtain
\[
 W(v) = \sum_{z=1}^{2Z-1} W_z(v) + W_0(v) = \sum_{z=1}^{2Z-1} V_z(v) + O (\Xi)
 + O \left( (\log N) \frac{N^{\gamma}}{Z} \right) ,
\]
where
\begin{equation*}
 \Xi  = \frac{vN^{\gamma - 1}}{Z} \sum_{P < p \leq 2P}(\log p)\sum_{z = 1}^{2Z - 1} g_{z}(p^{c}) .
\end{equation*}

\bigskip

Now we use \eqref{3}, \eqref{17.02.eq80} and Chebyshev's prime number theorem to find that
\[
 \Xi \ll \frac{vN^{ 2 \gamma - 1}}{Z}
\]
and therefore
\begin{equation}\label{85}
W(v) = \sum_{z = 1}^{2Z - 1}{V_z(v)} + O\left ( \frac{vN^{2\gamma - 1}}{Z}\right )
  + O\left ( (\log N)\frac{N^{\gamma }}{Z}\right ).
\end{equation}

\bigskip

From this point onwards we assume that
\begin{equation}\label{25}
 v = \frac{h}{d} , \qquad \text{where} \qquad 1 \le d \le D , \quad 1 \le h \le H .
\end{equation}
Then using \eqref{28} we see that
$ v N^{2\gamma - 1} \ll  N^{\gamma } (\log N)^3 $, hence 
formula \eqref{85} can be written as
\begin{equation} \label{18.02.eq10}
 W(v) = \sum_{z = 1}^{2Z - 1}{V_z(v)} + O\left ( (\log N)^3\frac{N^{\gamma }}{Z}\right ).
\end{equation}
From (\ref{28}), (\ref{83}) and \eqref{18.02.eq10} we find
\begin{equation}\label{95}
\Sigma _{j} \ll \frac{N^{2\gamma - 1}}{(\log N)^{2}} + \sum_{d \leq D}\sum_{h \leq H}{\frac{1}{h}}
\sum_{z = 1}^{2Z - 1}{|V_z(v)|} + \sum_{d \leq D}\sum_{h \leq H}{\frac{1}{h}{\frac{N^{\gamma }}{Z}}(\log N)^{3}}.
\end{equation}
We choose $Z$ such that
\begin{equation}\label{38}
  Z \asymp dN^{1 - \gamma }(\log N)^{7}.
\end{equation}
From \eqref{6} and \eqref{10050} it follows that the condition \eqref{86} holds. 
Consequently from \eqref{95} and \eqref{38} we find
\begin{equation}\label{89}
\Sigma_{j} \ll \frac{N^{2\gamma - 1}}{(\log N)^{2}} + 
\sum_{d \leq D}\sum_{h \leq H}{\frac{1}{h}}\sum_{z = 1}^{2Z - 1}{|V_z(v)|}.
\end{equation}

\bigskip

Now, we consider the sum $V_z(v)$ defined in \eqref{35}, where $v$ satisfies \eqref{25}.
By \eqref{17.02.eq60} we find that
\begin{align*}
V_z(v)
 & =
  \sum_{P < p \leq 2P}{(\log p)
  \left ( \sum_{|r| \leq Z (\log N)^{4}}{\beta _{r}^{(z)} \, e(rp^{c})}\right )
  e\left ( v\left ( N + j - p^{c} + \frac{z}{2Z}\right )^{\gamma }\right )} + O(N^{-10}) 
  \\
  & \\
 & =
  \sum_{|r| \leq Z (\log N)^{4}}{\beta _{r}^{(z)}} \,
  U \left(  N + j + \frac{z}{2Z} , r, v \right)
    + O(N^{-10}),
\end{align*}
where
\begin{equation}\label{39}
  U = U(T, r, v) = \sum_{P < p \leq 2P}{(\log p) \, e(rp^{c} + v(T - p^{c})^{\gamma })}.
\end{equation}

\bigskip

We would like to point out that when $0 \le j \le 1$ and $1 \le z \le 2Z-1$ 
then
\[
  N + j + \frac{z}{2Z} \in [N, N+2] .
\]
Furthermore, when we also take into account \eqref{17.02.eq70} and \eqref{38} we obtain
\[
  V_z(v) \ll N^{-10} + \frac{1}{Z}  \sum_{|r| \leq R }
  \; \sup_{T \in [N, N+2]} | U(T, r, v) | .
\]
where
\begin{equation} \label{40}
  R = d N^{1 - \gamma } (\log N)^{12}
\end{equation}
When we substitute the above expression for $V_z(v)$ in formula \eqref{89} we find that
\begin{equation}\label{98}
|\Sigma_{1}| + |\Sigma_{2}|
  \ll \frac{N^{2\gamma - 1}}{(\log N)^{2}} + \sum_{d \leq D} \; \sum_{h \leq H} {\frac{1}{h}} \;
\sum_{|r| \leq R} \; {\sup_{T \in [N, N+2]}{|U(T, r, v)|}} .
\end{equation}

\subsection{Application of Vaughan's identity}

Let us introduce the notations
\begin{equation}\label{41}
\phi (t) = rt^{c} + v(T - t^{c})^{\gamma },
\end{equation}
\begin{equation}
f(m, l) = \phi (ml) = r(ml)^{c} + v(T - (ml)^{c})^{\gamma }. \label{42}
\end{equation}
For the sum $U$, defined in \eqref{39}, we have
\[
 U = \sum_{P < n \le 2 P } \Lambda(n) \, e (\phi(n)) + O\left ( P^{\frac{1}{2}}\right ) .
\]
Now we apply Vaughan's identity (see \cite{5}) and find that
\begin{equation}\label{99}
U = U_{1} - U_{2} - U_{3} - U_{4} + O\left ( P^{\frac{1}{2}}\right ),
\end{equation}
where
\begin{align}
U_{1}
  & =
  \sum_{m \leq P^{\frac{1}{3}}}{\mu (m)}\sum_{\frac{P}{m} < l \leq \frac{2P}{m}}
  {(\log l) \, e(f(m, l))}, 
  \label{43} \\
  & \notag \\
U_{2}
  & =
  \sum_{m \leq P^{\frac{1}{3}}}{c(m)}\sum_{\frac{P}{m} < l \leq \frac{2P}{m}}{e(f(m, l))}, 
  \label{44} \\
  & \notag \\
U_{3}
  & =
  \sum_{P^{\frac{1}{3}} < m \leq P^{\frac{2}{3}}}{c(m)}\sum_{\frac{P}{m} < l \leq \frac{2P}{m}}{e(f(m, l))}, 
  \label{45} \\
  & \notag \\
U_{4}
  & =
  \mathop{\sum \sum}_{\substack{P < ml \leq 2P \\ m > P^{\frac{1}{3}}, \; \; l > P^{\frac{1}{3}}}}
  {a(m) \, \Lambda (l) \, e(f(m, l))}, 
   \label{46}
\end{align}
and where
\begin{equation}\label{47}
|c(m)| \leq \log m \hspace{5mm} \text{ and } \hspace{5mm} |a(m)| \leq \tau (m).
\end{equation}
Hence from \eqref{6}, \eqref{10050}, \eqref{28}, \eqref{40}, \eqref{98} and \eqref{99} we have
\begin{equation}\label{91}
|\Sigma_{1}| + |\Sigma_{2}| \ll \frac{N^{2\gamma - 1}}{(\log N)^{2}} + \sum_{i = 1}^{4}{\Omega _{i}},
\end{equation}
where
\begin{equation}\label{49}
\Omega_{i} = \sum_{d \leq D} \; \sum_{h \leq H}{\frac{1}{h}}
   \; \sum_{|r| \leq R} \; {\sup_{T \in [N, N + 2]}}{|U_{i}|}.
\end{equation}

By \eqref{91} and \eqref{49}, in order to prove that \eqref{10040} is satisfied, 
it suffices to show that 
\begin{equation}\label{90}
\Omega_{i} \ll \frac{N^{2\gamma - 1}}{(\log N)^{2}} 
  \qquad \text{for} \qquad  i = 1, 2, 3, 4 .
\end{equation}

\subsection{The estimation of the sums $\Omega_{1}$ and $\Omega_{2}$}

We begin with the study of $\Omega_{2}$. 
From \eqref{42} we get
\begin{equation}\label{51}
  f_{ll}''(m, l) = \pi_1 - \pi_2 ,
\end{equation}
where
\begin{equation}\label{10.03.2016.eq10}
  \pi_{1} = m^{2}rc(c-1)(ml)^{c-2} , \qquad
  \pi_{2} = m^{2}v(c-1)T(ml)^{c-2}(T-(ml)^{c})^{\gamma - 2} .
\end{equation}

\bigskip

From \eqref{3}, \eqref{25} and the conditions
\begin{equation}\label{10.03.2016.eq20}
 P < ml \leq 2P , \qquad N \leq T \leq N+2
\end{equation}
we have
\begin{equation}\label{52}
   |\pi_{1}| \asymp |r|m^{2}N^{1 - 2\gamma } \qquad \text {and} \qquad
    \pi_{2}  \asymp vm^{2}N^{-\gamma }.
\end{equation}

\bigskip

It follows from \eqref{51} and \eqref{52} that
there exists sufficiently small constant $\alpha _{0} > 0$ such that if
$|r| \leq \alpha _{0}vN^{\gamma - 1}$, then $|f_{ll}''| \asymp vm^{2}N^{-\gamma }$.

\bigskip

Similarly, from \eqref{51} and \eqref{52} we conclude that
there exists sufficiently large constant $A_{0} > 0$ such that if
$|r| \geq A_{0}vN^{\gamma - 1}$, then $|f_{ll}''| \asymp |r|m^{2}N^{1 - 2\gamma }$.

\bigskip

We divide the sum $\Omega_{2}$ into four sums according to the value of $r$ as follows:
\begin{equation}\label{129}
\Omega _{2} = \Omega _{2, 1} + \Omega _{2, 2} + \Omega _{2, 3} + \Omega _{2, 4},
\end{equation}
where
\begin{align}
   & \text{in  } \; \Omega _{2, 1} : \qquad |r| \leq \alpha_{0}vN^{\gamma - 1},
       \label{10.03.2016.eq100} \\
       & \notag \\
   & \text{in  } \;  \Omega _{2, 2} :  \qquad - A_{0}vN^{\gamma - 1} < r < - \alpha_{0}vN^{\gamma - 1},
      \label{10.03.2016.eq110} \\
      & \notag \\
   & \text{in  } \;  \Omega _{2, 3} :  \qquad \alpha_{0}vN^{\gamma - 1} < r < A_{0}vN^{\gamma - 1},
      \label{10.03.2016.eq120} \\
      & \notag \\
   & \text{in  } \;  \Omega _{2, 4} :  \qquad A_{0}vN^{\gamma - 1} \le |r| \le R.
      \label{10.03.2016.eq130}
\end{align}
We note that from \eqref{28}, \eqref{25} and \eqref{40} it follows that
\begin{equation} \label{10.04.2016.eq10}
 v N^{\gamma - 1} \ll  (\log N)^3 \ll  \frac{R}{\log N} .
\end{equation}

\bigskip

Let us consider $\Omega _{2, 4}$ first. We have
\begin{equation}\label{11.03.2016.eq20}
\Omega_{2, 4} = \sum_{d \leq D} \; \sum_{h \leq H} {\frac{1}{h}}
   \; \sum_{A_{0}vN^{\gamma - 1} \le |r| \le R } \; {\sup_{T \in [N, N + 2]}}{|U_{2}|}.
\end{equation}
We shall estimate the sum $U_2$, defined by \eqref{44}, 
provided that the condition \eqref{10.03.2016.eq130} holds.
We recall that the constant $A_0$ is chosen is such a way,
that if $|r| \geq A_{0}vN^{\gamma - 1}$, then uniformly for
$l \in \left( \frac{P}{m}, \frac{P}{m} \right]$ we have
$|f_{ll}''(m, l)| \asymp |r|m^{2}N^{1 - 2\gamma }$.
Hence we are in a position to use again Van der Corput's theorem (see \cite[Chapter 1, Theorem~5]{4})
and having also in mind \eqref{3} we obtain
\[
\sum_{\frac{P}{m}< l \leq \frac{2P}{m}}{e(f(m, l))}
   \ll
        \frac{P}{m}(|r|m^2N^{1-2\gamma })^{\frac{1}{2}} + (|r|m^{2}N^{1-2\gamma })^{-\frac{1}{2}}  \ll
        |r|^{\frac{1}{2}} \, N^{\frac{1}{2}} .
\]
Then from \eqref{44} and \eqref{47} we find
\begin{equation} \label{100}
   U_{2} \ll 
          |r|^{\frac{1}{2}} \, N^{\frac{1}{2} +\frac{\gamma }{3}} \, (\log N) .
\end{equation}
We substitute this expression for $U_2$ in \eqref{11.03.2016.eq20} and use
\eqref{28}, \eqref{25} and \eqref{40} to get
\[
\Omega_{2, 4}
     \ll
       D^{\frac{5}{2}} \, N^{2-\frac{7\gamma }{6} + \varepsilon} .
\]
Hence from \eqref{6} and \eqref{10050} we obtain
\begin{equation}\label{101}
\Omega _{2, 4} \ll \frac{N^{2\gamma - 1}}{(\log N)^{2}}.
\end{equation}

\bigskip

We carry on with $\Omega _{2, 3}$.
We have to study the sum $U_2$ defined by \eqref{44} provided that the condition
\eqref{10.03.2016.eq120} holds.
To do so we use \eqref{42} and compute
\begin{align}\label{54}
f_{lll}'''(m, l)
    & =
       m^{3}rc(c-1)(c-2)(ml)^{c-3}  \notag \\
       & \notag \\
    & \qquad +
       m^{3}v(c-1)T(T - (ml)^{c})^{\gamma - 3}(ml)^{c-3}((2-c)T - (c+1)(ml)^{c}).
\end{align}

\bigskip

From \eqref{51}, \eqref{10.03.2016.eq10} and \eqref{54}
we find
\[
    (c-2)f_{ll}''(m, l) - lf_{lll}'''(m, l)
    = v(c-1)(2c-1) T m^{2c} l^{2c - 2} \left(T - (ml)^c \right)^{\gamma - 3} .
\]
The above, together with \eqref{3} and \eqref{10.03.2016.eq20} implies
\[
  \left| (c-2)f_{ll}''(m, l) - lf_{lll}'''(m, l) \right| \asymp v m^{2}N^{-\gamma }.
\]
Therefore, there exists  $\kappa _{0} > 0$,
depending only on the constant $c$, such that for every $l \in \left ( \frac{P}{m}, \frac{2P}{m} \right ]$
at least one of the following inequalities holds:
\begin{equation} \label{105}
|f_{ll}''(m, l)| \geq \kappa _{0} vm^{2}N^{-\gamma },
\end{equation}
or
\begin{equation}\label{106}
|f_{lll}'''(m, l)| \geq \kappa _{0} vm^{3}N^{-2\gamma }.
\end{equation}

\bigskip

We are going to show that
the interval $ \left ( \frac{P}{m}, \frac{2P}{m} \right ]$ can be divided into
at most $7$ intervals such that if $J$ is one of them,
then 
at least one of the following statements holds:
\begin{equation} \label{31.03.2016.eq10}
 \text{We have \eqref{105} for all} \;\; l \in J .
\end{equation}

\begin{equation}\label{31.03.2016.eq20}
 \text{We have \eqref{106} for all} \;\; l \in J .
\end{equation}

\bigskip

To establish this it is enough to show that the equation
$|f_{ll}''(m, l)| = \kappa _{0} vm^{2}N^{-\gamma } $ has at most $6$ 
solutions in real numbers
$l \in \left ( \frac{P}{m}, \frac{2P}{m} \right )$.
Hence, it is enough to show that if $C$ does not depend on $l$ then
the equation $ f_{ll}''(m, l) = C $ has at most $3$ solutions in real numbers
$l \in \left ( \frac{P}{m}, \frac{2P}{m} \right )$.
According to Rolle's theorem, between any two solutions of the last equation
there is a solution (in real numbers $l$) of the equation $f_{lll}''' (m, l) = 0$.
So, we use \eqref{54} to conclude that it is enough to show that
\[
   vT(T-(ml)^{c})^{\gamma - 3}((2-c)T - (c+1)(ml)^{c}) =  rc(2-c)
\]
has at most $2$ solutions in $l \in \left ( \frac{P}{m}, \frac{2P}{m} \right ) $, which is equivalent
to the assertion that the equation
\[
   (T-X)^{\gamma - 3}((2-c)T - (c+1)X)= \frac{ rc(2-c) }{vT}
\]
has at most $2$ solutions in $X \in \left( P^c , (2P)^c \right)$.
Alternatively, instead of the last equation one can look at
\begin{equation} \label{10.03.2016.eq30}
  (\gamma - 3 ) \log (T - X) + \log \left( (2-c)T - (c+1)X \right) = \log \frac{ rc(2-c) }{vT} .
\end{equation}
Let $H(X)$ denote the function on the left side of \eqref{10.03.2016.eq30}.
From Rolle's theorem we know that between any two solutions of \eqref{10.03.2016.eq30}
there is a solution of $H'(X) = 0 $.
Since
\[
  H'(X) = \frac{3 - \gamma}{T-X} - \frac{c+1}{(2 - c )T - (c+1) X}
\]
it is easy to see that $H'(X)$ vanishes for at most $1$ value of $X$.
Therefore, \eqref{10.03.2016.eq30} has at most $2$ solutions in $X$
and our assertion is proved.

\bigskip

On the other hand, from \eqref{51}, \eqref{52} and \eqref{54}
we see that under the condition on $r$ imposed in \eqref{10.03.2016.eq120} we have
\[
 f_{ll}'' (m, l) \ll v m^2 N^{- \gamma} \qquad \text{and} \qquad
  f_{lll}''' (m, l) \ll v m^3 N^{- 2 \gamma} .
\]
Hence, the interval $ \left ( \frac{P}{m}, \frac{2P}{m} \right ]$ can be divided into
at most $7$ intervals such that if $J$ is one of them, then
at least one of the following assertions holds:
\begin{equation} \label{10.03.2016.eq160}
|f_{ll}''(m, l)| \asymp vm^{2}N^{-\gamma } \qquad \; \text{uniformly for} \qquad l \in J ,
\end{equation}

\begin{equation} \label{10.03.2016.eq170}
|f_{lll}'''(m, l)| \asymp  vm^{3}N^{-2\gamma } \qquad \text{uniformly for} \qquad l \in J .
\end{equation}

\bigskip

If \eqref{10.03.2016.eq160} is fulfilled, 
then we use Van der Corput's theorem (see \cite[Chapter~1, Theorem~5]{4})
for the second derivative and find
\begin{equation} \label{11.03.2016.eq65}
  \sum_{l \in J} e (f (m, l))
  \ll  \frac{P}{m} \left( v m^2 N^{-\gamma} \right)^\frac{1}{2}
      + \left( v m^2 N^{-\gamma} \right)^{-\frac{1}{2}}
  \ll  v^{\frac{1}{2}}N^{\frac{\gamma}{2}} + v^{-\frac{1}{2}}m^{-1}N^{\frac{\gamma}{2}} .
\end{equation}
In the case \eqref{10.03.2016.eq170} we apply Van der Corput's theorem for the third derivative
to get
\[
  \sum_{l \in J} e (f (m, l))
  \ll
   \frac{P}{m} \left( vm^{3}N^{-2\gamma } \right)^{\frac{1}{6}} +
    \left( \frac{P}{m} \right)^{\frac{1}{2}} \left( v m^{3} N^{-2\gamma } \right)^{-\frac{1}{6}}
  \ll v^{\frac{1}{6}}m^{-\frac{1}{2}}N^{\frac{2\gamma }{3}}
      + v^{-\frac{1}{6}}m^{-1}N^{\frac{5\gamma }{6}} .
\]
Hence, in either case $\sum_{l \in J} e (f (m, l)) $ can be estimated by the
sum of the expressions on the right sides of the inequalities above.
Therefore,
\begin{equation} \label{11.03.2016.eq75}
\sum_{\frac{P}{m} < l \leq \frac{2P}{m}}{e(f(m, l))}
    \ll    v^{\frac{1}{2}}N^{\frac{\gamma}{2}}
        + v^{-\frac{1}{2}}m^{-1}N^{\frac{\gamma}{2}}
        + v^{\frac{1}{6}}m^{-\frac{1}{2}}N^{\frac{2\gamma }{3}}
        + v^{-\frac{1}{6}}m^{-1}N^{\frac{5\gamma }{6}} .
\end{equation}
Then from \eqref{3}, \eqref{44} and \eqref{47} we find that
\begin{equation}\label{57}
 U_{2}
     \ll (\log N)^2 \left ( v^{\frac{1}{2}} N^{\frac{5\gamma}{6}}  + v^{-\frac{1}{2}} N^{\frac{\gamma }{2}}+
        v^{\frac{1}{6}}N^{\frac{5\gamma }{6}} + v^{-\frac{1}{6}}N^{\frac{5\gamma }{6}}\right ).
\end{equation}

\bigskip

We use \eqref{28}, \eqref{25}, \eqref{10.03.2016.eq120}, \eqref{10.04.2016.eq10}
and \eqref{57} to get
\begin{align} 
\Omega_{2, 3} 
   & = 
   \sum_{d \leq D} \; \sum_{h \leq H}{\frac{1}{h}}
   \; \sum_{\alpha_{0}vN^{\gamma - 1} < r < A_{0}vN^{\gamma - 1}} \; {\sup_{T \in [N, N + 2]}}{|U_{2}|}
    \notag \\
    & \notag \\
   & \ll
   N^{\varepsilon}
   \sum_{d\leq D} \; \sum_{h \leq H}{\frac{1}{h}}
            \left ( \frac{ h^{\frac{1}{2}} N^{\frac{5\gamma}{6}} }{d^{\frac{1}{2}}} +
            \frac{ h^{-\frac{1}{2}} N^{\frac{\gamma }{2}}  }{d^{-\frac{1}{2}}}
            + \frac{h^{\frac{1}{6}}N^{\frac{5\gamma }{6}}}{d^{\frac{1}{6}}} +
            \frac{h^{-\frac{1}{6}}N^{\frac{5\gamma }{6}}}{d^{-\frac{1}{6}}}\right ) 
            \notag \\
            & \notag \\
    & \ll  
    N^{\varepsilon}
            \left ( DN^{\frac{1}{2} + \frac{\gamma }{3}} + D^{\frac{3}{2}}N^{\frac{\gamma }{2}}
            + DN^{\frac{1}{6} + \frac{2\gamma }{3}} + D^{\frac{7}{6}}N^{\frac{5\gamma }{6}} \right )
\end{align}
and from \eqref{6} and \eqref{10050} we deduce that
\begin{equation}\label{102}
\Omega _{2, 3} \ll \frac{N^{2\gamma - 1}}{(\log N)^{2}}.
\end{equation}

\bigskip

Let us consider $\Omega _{2, 1}$.
We have chosen the constant $\alpha_0$ in such a way, that from \eqref{10.03.2016.eq20}
and from the condition on $r$ imposed in \eqref{10.03.2016.eq100}
it follows that
$| f_{ll}''(m, l) |\asymp vm^{2}N^{-\gamma }$
uniformly for $l \in \left( \frac{P}{m} , \frac{2P}{m} \right]$.
Hence the sum $\sum_{\frac{P}{m} < l \le \frac{2P}{m} } e \left( f(m, l) \right)$
can be estimated by the expression on the right side of
\eqref{11.03.2016.eq65} and certainly the estimate \eqref{11.03.2016.eq75} holds again.
From this observation we see that $\Omega _{2, 1}$ can be estimated
in the same way as $\Omega_{2, 3}$, i.e.
\begin{equation}\label{130}
\Omega _{2, 1} \ll \frac{N^{2\gamma - 1}}{(\log N)^{2}}.
\end{equation}

\bigskip

The sum $\Omega _{2, 2}$ can be studied in the same way.
From \eqref{51} -- \eqref{10.03.2016.eq20}
and \eqref{10.03.2016.eq110} it follows that 
$|f''_{ll}(m, l)| \asymp v m^2 N^{- \gamma}$ and hence the estimate \eqref{11.03.2016.eq75}
is correct again. Therefore
\begin{equation}\label{131}
\Omega _{2, 2} \ll \frac{N^{2\gamma - 1}}{(\log N)^{2}}.
\end{equation}

\bigskip

From \eqref{129}, \eqref{101} and \eqref{102} -- \eqref{131} we conclude that
\begin{equation}\label{136}
\Omega _{2} \ll \frac{N^{2\gamma - 1}}{(\log N)^{2}}.
\end{equation}

\bigskip

Consider now $\Omega_{1}$.
For $U_{1}$ defined by \eqref{43}, we use Abel's transformation to get rid of the factor $\log l$
in the inner sum. Then we proceed as in the estimation of $\Omega_{2}$
to obtain
\begin{gather}\label{18.02.eq50}
\Omega _{1} \ll \frac{N^{2\gamma - 1}}{(\log N)^{2}}.
\end{gather}

\subsection{The estimation of the sums $\Omega_{3}$ and $\Omega_{4}$ and the end of the proof}

Consider the sum $\Omega_4$, defined in \eqref{49}.
We divide the sum $U_{4}$ given by \eqref{46} into  $O \left(\log N \right) $ sums of the form
\begin{equation}\label{60}
W_{M,L} = \sum_{L < l \leq 2L}{b(l)}
        \sum_{\substack {M < m \leq 2M \\ \frac{P}{l} < m \leq \frac{2P}{l}}}{a(m) e(f(m,l))},
\end{equation}
where
\begin{equation}\label{61}
a(m) \ll N^{\varepsilon }\text{, }\hspace{5mm}
b(l) \ll N^{\varepsilon }\text{, }\hspace{5mm}
P^{\frac{1}{3}} \leq M \leq P^{\frac{1}{2}} \ll L \ll P^{\frac{2}{3}}\text{, }\hspace{5mm}
ML \asymp P.
\end{equation}

\bigskip

From \eqref{60}, \eqref{61} and Cauchy's inequality we find that
\begin{equation}\label{62}
|W_{M,L}|^{2} \ll N^{\varepsilon }L\sum_{L < l \leq 2L}
                \left | \sum_{M_{1} < m \leq M_{2}}{a(m) \, e(f(m, l))}\right |^{2},
\end{equation}
where
\begin{equation}\label{63}
M_{1} = \max{\left ( M, \frac{P}{l}\right )}, \hspace{10mm} M_{2} = \min{\left ( 2M, \frac{2P}{l} \right )}.
\end{equation}
Now we apply the well-known inequality
\begin{equation}\label{93}
\left | \sum_{a < m \leq b}{\xi (m)}\right |^{2} \leq \frac{b-a+Q}{Q}\sum_{|q|<Q}\left ( 1 - \frac{|q|}{Q}\right )
\sum_{\substack {m \in (a, b] \\ m \in (a-q, b-q]}}{\xi (m+q)} \, \overline{\xi (m)},
\end{equation}
where $Q \in \mathbb{N}$, $a, b \in \mathbb{R}$, $1 \leq b-a$ and $\xi (m)$ is any complex function.
(A proof can be found in~\cite[Lemma 8.17]{9}).
In our setting $\xi(m) = a(m) \, e (f(m, l))$,
$a = M_1, b = M_2$.
The exact value of $Q$ will be chosen later. For now we only require that
\begin{equation} \label{12.03.2016.eq10}
  Q \le M .
\end{equation}
Then we find
\begin{align*}
|W_{M,L}|^{2}
        & \ll
             N^{\varepsilon }L\sum_{L < l \leq 2L}{\frac{M}{Q}}\sum_{|q|\leq Q}
             {\left ( 1 - \frac{|q|}{Q}\right )}
              \\
              & \\
        & \qquad \times
             \sum_{\substack {M_{1} < m \leq M_{2} \\ M_{1}< m + q \leq M_{2}}}
                                    {a(m+q) \, \overline{a(m)} \, e(f(m+q, l)-f(m, l))} .
\end{align*}
We estimate the contribution coming from the terms with $q=0$, then we change the order of summation 
and using \eqref{61} and \eqref{63} we find
\begin{equation} \label{64}
|W_{M,L}|^{2}
         \ll
             \frac{N^{\varepsilon }(LM)^{2}}{Q} + \frac{N^{\varepsilon }LM}{Q}\sum_{0< |q| \leq Q}
             \sum_{\substack {M < m \leq 2M \\ M < m+q \leq 2M}}
             \left|
             {\sum_{L_{1} < l \leq L_{2}}
             e \left(  Y_{m, q}(l) \right) } \right| ,
\end{equation}
where
\begin{equation}\label{65}
L_{1} = \max \left ( L, \frac{P}{m}, \frac{P}{m+q}\right ), \hspace{10mm}
L_{2} = \min \left ( 2L, \frac{2P}{m}, \frac{2P}{m+q}\right )
\end{equation}
and
\begin{equation}\label{66}
Y(l) = Y_{m, q}(l) = f(m+q, l) - f(m, l).
\end{equation}

\bigskip

It is now easy to see that the sum over negative $q$ in formula \eqref{64} is equal to the sum over 
positive $q$, hence we obtain
\begin{equation} \label{31.03.2016.eq75}
|W_{M,L}|^{2}
        \ll
             \frac{N^{\varepsilon }(LM)^{2}}{Q} + \frac{N^{\varepsilon }LM}{Q}\sum_{1 \le q \leq Q}
             \;
             \sum_{M < m \leq 2M - q } \,
             \left|
             {\sum_{L_{1} < l \leq L_{2}}
             e \left(  Y_{m, q}(l) \right) } \right| .
\end{equation}

\bigskip

Consider the function $Y(l)$.
Using \eqref{41}, \eqref{42} and \eqref{66}
we find that
\begin{equation}\label{67}
Y(l) = \int_{m}^{m+q}{f_{t}'(t, l)dt}=\int_{m}^{m+q}{l\phi '(tl)dt} 
\end{equation}
and therefore
\begin{equation}\label{68}
Y''(l) = \int_{m}^{m+q}{(2t\phi '' (tl) + lt^{2}\phi ''' (tl))dt}, \qquad
Y'''(l) = \int_{m}^{m+q}{(3t^{2}\phi ''' (tl) + lt^{3}\phi ^{(4)}(tl))dt}.
\end{equation}
From \eqref{41} and \eqref{68} we get
\begin{align}\label{70}
Y''(l) = \int_{m}^{m+q}{(\Phi_{1}(t) - {\Phi_2}(t))dt},
\end{align}
where
\begin{align}
 \Phi_{1} (t)
  & = 
  rc^{2}(c-1)t^{c-1}l^{c-2} ,
   \label{01.04.2016.eq20} \\
   & \notag \\
 \Phi_{2} (t)
  & = 
  v(c-1)Tt^{c-1}l^{c-2} \left( T-(tl)^{c} \right)^{\gamma -3} \left( cT + (c-1) (tl)^c \right).
  \label{01.04.2016.eq30}
\end{align}

If $t \in [m, m+q]$ then $tl \asymp P$. 
Thus, by \eqref{3} and the condition
\begin{equation} \label{01.04.2016.eq10}
  N \le T \le N+2
\end{equation}
we find that uniformly for $t \in [m, m+q]$ we have
\begin{equation}\label{71}
    |\Phi_{1} (t) |\asymp |r| \, m \, N^{1-2\gamma } \qquad \text {and} \qquad
    \Phi_{2} (t) \asymp v \, m \, N^{-\gamma }.
\end{equation}

\bigskip

From \eqref{70} and \eqref{71} we see that there exists a sufficiently small constant 
$\alpha _{1} > 0$ such that if
$|r| \leq \alpha _{1}vN^{\gamma - 1}$, then $|Y''(l)| \asymp qvmN^{-\gamma }$.
Similarly, we conclude that there exists a sufficiently large constant $A_{1} > 0$
such that if $|r| \geq A_{1}vN^{\gamma - 1}$, then $|Y''(l)| \asymp |r|qmN^{1 - 2\gamma }$.
Hence, it makes sense to divide the sum $\Omega_{4} $ 
into four sums according to the value of $r$ as follows:
\begin{equation}\label{132}
\Omega_{4}  = \Omega_{4, 1} + \Omega_{4, 2} + \Omega_{4, 3} + \Omega_{4, 4},
\end{equation}
where
\begin{align}
   & \text{in  } \; \Omega _{4, 1} : \qquad |r| \leq \alpha_{1}vN^{\gamma - 1},
       \label{15.03.2016.eq10} \\
       & \notag \\
   & \text{in  } \;  \Omega _{4, 2} :  \qquad - A_{1}vN^{\gamma - 1} < r < - \alpha_{1}vN^{\gamma - 1},
      \label{15.03.2016.eq11} \\
      & \notag \\
   & \text{in  } \;  \Omega _{4, 3} :  \qquad \alpha_{1}vN^{\gamma - 1} < r < A_{1}vN^{\gamma - 1},
      \label{15.03.2016.eq12} \\
      & \notag \\
   & \text{in  } \;  \Omega _{4, 4} :  \qquad A_{1}vN^{\gamma - 1} \le |r| \le R.
      \label{15.03.2016.eq13}
\end{align}

\bigskip

Let us consider $\Omega _{4, 4}$ first.
From \eqref{49} and \eqref{15.03.2016.eq13} we have
\begin{equation} \label{15.03.2016.eq144}
\Omega_{4, 4}
        \ll (\log N)\sum_{d \leq D} \; \sum_{h \leq H}{\frac{1}{h}}
        \; \sum_{A_{1}vN^{\gamma - 1} \le |r| \le R} \; 
        {\sup_{\substack{T \in [N, N + 2] \\ M, L \; : \; \eqref{61} }}}{|W_{M, L}|}. 
\end{equation}
(The supremum is taken over $T \in [N, N+2]$ and $M$, $L$ satisfying the conditions imposed in \eqref{61}).

\bigskip

Consider the sum $W_{M, L}$. We already mentioned that if
$|r| \geq A_{1}vN^{\gamma - 1}$, then uniformly for 
$l \in \left ( L_1, L_2 \right ]$ we have
$Y''(l) \asymp |r| \, q m N^{1 - 2\gamma }$. Hence we can use Van der Corput's theorem 
(see \cite[Chapter~1, Theorem~5]{4}) for the second derivative and
by \eqref{3}, \eqref{61} and \eqref{65} we obtain
\[
\sum_{L_{1} < l \leq L_{2}}{e(Y(l))}
        \ll
             L (|r|qmN^{1-2\gamma })^{\frac{1}{2}} + (|r|qmN^{1-2\gamma })^{-\frac{1}{2}}\\
        \ll
            |r|^{\frac{1}{2}} q^{\frac{1}{2}} M^{-\frac{1}{2}} N^{\frac{1}{2}} .
\]
Then from \eqref{3}, \eqref{61} and \eqref{31.03.2016.eq75} we find
\begin{equation} \label{75}
 W_{M, L}
        \ll
            N^{\varepsilon } \left ( {N^{\gamma }} Q^{-\frac{1}{2}}
            + |r|^{\frac{1}{4}} \, Q^{\frac{1}{4}} \, N^{\frac{1}{4} + \frac{5 \gamma}{8}} \right) .
\end{equation}
From \eqref{28}, \eqref{40} and \eqref{15.03.2016.eq144} we have
\begin{align*}
\Omega _{4, 4}
      & \ll
            N^{\varepsilon } 
            \sum_{d \leq D} \sum_{h \leq H} {\frac{1}{h}}
                 \sum_{|r| \leq dN^{1-\gamma } (\log N)^{12}} 
             \left ( N^{\gamma } Q^{-\frac{1}{2}} 
            + |r|^{\frac{1}{4}} \, Q^{\frac{1}{4}} \, N^{\frac{1}{4} + \frac{5\gamma }{8}} \right)
             \\
             & \\
      & \ll
            N^{\varepsilon }\left ( D^{2} N Q^{-\frac{1}{2}} + Q^{\frac{1}{4}}D^{\frac{9}{4}}
            N^{\frac{3}{2}-\frac{5\gamma }{8}} 
             \right ).
\end{align*}
We choose 
\begin{equation} \label{02.04.2016.eq10}
   Q = \left[ D^{-\frac{1}{3}}N^{\frac{5\gamma }{6} - \frac{2}{3}} \right] . 
\end{equation}
It is now easy to verify that the condition \eqref{12.03.2016.eq10} holds.
Hence, from \eqref{6} and \eqref{10050} we obtain
\begin{equation}\label{103}
\Omega _{4, 4} \ll \frac{N^{2\gamma - 1}}{(\log N)^{2}}.
\end{equation}

\bigskip

Let us now consider  $\Omega_{4, 3}$. From \eqref{49} and \eqref{15.03.2016.eq12} we have
\begin{equation}
\Omega_{4, 3}
       \ll (\log N)\sum_{d \leq D} \; \sum_{h \leq H}{\frac{1}{h}}
        \; \sum_{\alpha_{1}vN^{\gamma - 1} < r < A_{1}vN^{\gamma - 1}} \; 
        {\sup_{\substack{T \in [N, N + 2] \\ M, L \; : \;  \eqref{61} }}}{|W_{M, L}|}. 
        \label{15.03.2016.eq14}
\end{equation}
Consider the sum $W_{M, L}$ from the expression in the above formula. 
Using \eqref{70} -- \eqref{01.04.2016.eq30} we find
\begin{equation} \label{01.04.2016.eq45}
  Y'''(l) = \int_m^{m+q} ( \Psi_1 (t) + \Psi_2 (t) ) \, d t ,
\end{equation}
where
\begin{align}
  \Psi_1 (t)
   & =
    r \, c^2 \, (c-1) \, (c-2) \, t^{c-1} \, l^{c-3} ,
    \label{01.04.2016.eq50} \\
    & \notag \\
   \Psi_2 (t)
     & =  
     v \, (c-1) \, T \, t^{c-1} \, l^{c-3} \,
     \left( T - (tl)^c \right)^{\gamma - 4} 
      \notag \\
      & \notag \\
    & \qquad \qquad \times  
     \Big( c(2-c) \, T^2 + (-4 c^2 + 3 c - 2) \, T \, (tl)^c  + (1 - c^2) \, (tl)^{2c} \Big).
    \label{01.04.2016.eq55}
\end{align}
From 
\eqref{70} -- \eqref{01.04.2016.eq30} and
\eqref{01.04.2016.eq45} -- \eqref{01.04.2016.eq55}
we obtain
\[
  l Y'''(l) + (2-c) Y''(l)
  = - \int_m^{m+q} \Theta (t) \, d t ,  
\]
where
\[
  \Theta (t) = v(c-1) T t^{2c-1} l^{2c - 2 } 
  \left( T - (tl)^c \right)^{\gamma - 4} \,
  \left( 2c (2c-1) T + (2c^2 - 3c +1) (tl)^c \right) .
\]
Using \eqref{3}, \eqref{61} and \eqref{01.04.2016.eq10} 
we find that uniformly for $t \in [m, m+q]$ we have
\[
  \Theta(t) \asymp v m N^{- \gamma} .
\]
Hence
\[
  | (2-c)Y''(l) + l Y'''(l) | \asymp q m v N^{-\gamma }
\]
uniformly for $l \in \left( L_1 , L_2 \right]$.
Therefore, there exists $\kappa _{1} > 0$ which depends only on $c$ 
and such that, at least one of the following inequalities holds for every $l \in \left ( L_1, L_2 \right ]$:
\begin{align} 
  |Y''(l)| 
   & \geq 
   \kappa_{1}  v  q  m  N^{-\gamma },
  \label{73} \\
  & \notag \\
  |Y'''(l)| 
   & \geq 
   \kappa_{1}  v  q  m^{2} N^{-2\gamma }.
  \label{74}
\end{align}

\bigskip

The next step is to show that the interval
$\left ( L_1, L_2 \right ]$ can be divided into at most 13 intervals 
such that if $J$ is one of them, 
then at least one of the following assertions holds:
\begin{align} 
   &
   \text{We have \eqref{73} for all }  l \in J .
    \label{01.04.2016.eq95} \\
    & \notag \\
   & 
   \text{We have \eqref{74} for all }  l \in J .
    \label{01.04.2016.eq100}
\end{align}
To establish this, it suffices to show
that the equation $|Y''(l)| = \kappa _{1}vqmN^{-\gamma }$ has at most 12 solutions in real numbers
$l \in \left ( L_1, L_2 \right )$. Hence, it is enough to show that if $C$ does not depend on $l$, 
then the equation $Y''(l) = C$ has at most $6$ solutions in real numbers $l \in \left ( L_1, L_2 \right )$. 
According to Rolle's theorem, between any two solutions of the last equation 
there is a solution of the equation 
\begin{equation} \label{01.04.2016.eq120}
  Y'''(l) = 0 . 
\end{equation}
Hence, it is enough to show that \eqref{01.04.2016.eq120}
has at most $5$ solutions in real numbers $l \in (L_1, L_2)$.

\bigskip

By \eqref{54} and \eqref{66} one can easily see that
\eqref{01.04.2016.eq120} is equivalent to
the equation
\begin{align}
   &
  (m+q)^c \big( T - (m+q)^c l^c \big)^{\gamma - 3} 
  \big( (2 - c) T - (c+1) (m+q)^c l^c \big)
   \notag \\
   & \notag \\
  & \qquad 
  - 
  m^c \left( T - m^c l^c \right)^{\gamma - 3} 
  \big( (2 - c) T - (c+1) m^c l^c \big)
  = \frac{rc (2-c) \big( (m+q)^c - m^c \big) }{vT} .
  \notag
\end{align}

\bigskip

Let $X = l^c$. Define
\begin{align}
 \mathcal F (X)
  & = 
   (m+q)^c \big( T - (m+q)^c X \big)^{\gamma - 3} 
  \big( (2 - c) T - (c+1) (m+q)^c X \big)
    \notag \\
    & \notag \\
   & \qquad 
     - 
  m^c \left( T - m^c X \right)^{\gamma - 3} 
  \big( (2 - c) T - (c+1) m^c X \big) .
  \notag
\end{align}
It would be enough to show that if
$B$ does not depend on $X$, then the equation
$ \mathcal F(X) = B $
has at most $5$ solutions with $X \in \left( L_1^c, L_2^c \right)$.

\bigskip

Once more, we refer to Rolle's theorem to justify that 
it is enough to prove that the equation 
$ \mathcal F '(X)=0$ has no more than $4$ solutions with 
$X \in \left ( L_1^c , L_2^c \right )$.
One could write $ \mathcal F '(X)=0$ as
\begin{align*}
    & 
      (m+q)^{2c}\big( T-(m+q)^{c} X \big)^{\gamma - 4} \big( (4c+2\gamma -6)T + (2c-\gamma + 1)(m+q)^{c} X \big) 
       \\
       & \\
    &  \qquad = 
    m^{2c}\big(T-m^{c} X \big)^{\gamma - 4} \big( (4c+2\gamma -6)T + (2c-\gamma +1)m^{c} X \big) ,
\end{align*}
which, in turn, is equivalent to
\[
  \mathcal G(X) = \log m^c - \log (m+q)^c ,
\]
where
\begin{align}
  \mathcal G (X) 
   & = 
   (\gamma - 4) \log \big( T - (m+q)^c X \big) 
                      +  \log \big( (4 c + 2 \gamma - 6) T + (2 c + 1 - \gamma)  (m+q)^c X \big)
                      \notag \\
        & \notag \\
  & \qquad -
   (\gamma - 4) \log \big( T - m^c X \big) 
                      -  \log \big( (4 c + 2 \gamma - 6) T + (2 c + 1 - \gamma)  m^c X \big) .
    \notag                                         
\end{align}
By the same argument as before it is enough
to establish that
the equation 
\begin{equation} \label{01.04.2016.eq170}
 \mathcal G'(X) = 0 
\end{equation}
has at most $3$ solutions with $X \in \left ( L_1^c , L_2^c \right )$.
This can easily be shown because 
\begin{align}
  \mathcal G'(X) 
   & =
   \frac{(4 - \gamma) (m+q)^c }{T - (m+q)^c X} 
     + \frac{(2 c + 1 - \gamma) (m+q)^c}{(4c + 2 \gamma - 6) T + (2 c + 1 - \gamma) (m+q)^c X}
     \notag \\
     & \notag \\
   & \qquad -
     \frac{(4 - \gamma) m^c }{T - m^c X} 
     - \frac{(2 c + 1 - \gamma) m^c}{(4c + 2 \gamma - 6) T + (2 c + 1 - \gamma) m^c X} .
     \notag
\end{align}
Therefore, the number of solutions of \eqref{01.04.2016.eq170} 
does not exceed the number of roots of non-zero polynomial of degree at most $3$.

\bigskip

On the other hand, from \eqref{70}, \eqref{71}, \eqref{15.03.2016.eq12} and
\eqref{01.04.2016.eq45} -- \eqref{01.04.2016.eq55} we have
\begin{equation}\label{15.03.2016.eq15}
 Y'' (l) \ll vqmN^{- \gamma} \qquad \text{and} \qquad
 Y''' (l) \ll vqm^{2}N^{- 2 \gamma} .
\end{equation}
Hence, we conclude that the
interval $ \left ( L_1, L_2 \right ]$ can be divided into
at most $13$ intervals such that if $J$ is one of them, then
at least one of the following assertions holds:
\begin{align} 
  & 
  |Y''(l)| \asymp vqmN^{-\gamma } \qquad \;\;\; \text{uniformly for} \qquad l \in J ,
  \label{16.03.2016.eq160} \\
  & \notag \\
  &
  |Y'''(l)| \asymp vqm^{2}N^{-2\gamma } \qquad \text{uniformly for} \qquad l \in J .
\label{16.03.2016.eq170}
\end{align}

\bigskip

If  \eqref{16.03.2016.eq160} holds,
then we use \eqref{3}, \eqref{61} and Van der Corput's theorem (see \cite[Chapter~1, Theorem~5]{4})
for the second derivative  to get
\begin{align}
\sum_{l \in J}{e(Y(l))}
    & \ll 
    L (qvmN^{-\gamma })^{\frac{1}{2}} + (qvmN^{-\gamma })^{-\frac{1}{2}}
     \notag \\
    & \ll 
          q^{\frac{1}{2}} v^{\frac{1}{2}} L M^{\frac{1}{2}}N^{-\frac{\gamma }{2}}
        + q^{-\frac{1}{2}} v^{-\frac{1}{2}} M^{-\frac{1}{2}} N^{\frac{\gamma }{2}}. \label{16.03.2016.eq180}
\end{align}
In the case when \eqref{16.03.2016.eq170} is satisfied we
apply \eqref{3}, \eqref{61} and
Van der Corput's theorem for the third derivative to get
\begin{align}
\sum_{l \in J}{e(Y(l))}
    & \ll
         L (qvm^{2}N^{-2\gamma })^{\frac{1}{6}}
         + L^{\frac{1}{2}} (q v m^{2} N^{-2\gamma })^{-\frac{1}{6}} \notag \\
    & \ll
         q^{\frac{1}{6}} v^{\frac{1}{6}} L M^{\frac{1}{3}}N^{ - \frac{\gamma }{3}}
         + q^{-\frac{1}{6}} v^{-\frac{1}{6}} L^{\frac{1}{2} } M^{-\frac{1}{3}} N^{\frac{\gamma }{3}}. \label{16.03.2016.eq181}
\end{align}

\bigskip

Hence, in each case, $\sum_{l \in J}{e(Y(l))}$ can be estimated by the sum of the expressions on the right sides of the
inequalities \eqref{16.03.2016.eq180} and \eqref{16.03.2016.eq181}. Therefore, we obtain
\[
\sum_{L_1 < l \leq L_2} e(Y(l))
    \ll
     q^{\frac{1}{2}} v^{\frac{1}{2}} L M^{\frac{1}{2}}N^{-\frac{\gamma }{2}}
        + q^{-\frac{1}{2}} v^{-\frac{1}{2}} M^{-\frac{1}{2}} N^{\frac{\gamma }{2}}
         + q^{\frac{1}{6}} v^{\frac{1}{6}} L M^{\frac{1}{3}}N^{ - \frac{\gamma }{3}}
         + q^{-\frac{1}{6}} v^{-\frac{1}{6}} L^{\frac{1}{2} } M^{-\frac{1}{3}} N^{\frac{\gamma }{3}} .
\]
We use \eqref{3}, \eqref{61}, \eqref{65} and \eqref{31.03.2016.eq75}
and find that
\[
 W_{M, L}  \ll N^{\varepsilon }
    \left ( N^{\gamma } Q^{-\frac{1}{2}}
    + v^{\frac{1}{4}} Q^{\frac{1}{4}} N^{\frac{7\gamma }{8}}
    +  v^{-\frac{1}{4}} Q^{-\frac{1}{4}} N^{\frac{7\gamma }{8}}
    + v^{\frac{1}{12}} Q^{\frac{1}{12}} N^{\frac{11\gamma }{12}}
    + v^{-\frac{1}{12}} Q^{-\frac{1}{12}} N^{\frac{23\gamma }{24}}
    \right ).
\]

\bigskip

We apply the above estimate for $W_{M, L}$  in \eqref{15.03.2016.eq14}. 
Then, by \eqref{28} and \eqref{25} we obtain
\begin{align*}
\Omega _{4, 3}
    & \ll
         N^{\varepsilon } \; \sum_{d \leq D} \; \sum_{h \leq H}{\frac{1}{h}}
         \sum_{ r < A_{1} \log^3 N} \Big( N^{\gamma } Q^{-\frac{1}{2}}
         + \Big( \, \frac{h}{d} \, \Big)^{\frac{1}{4}}Q^{\frac{1}{4}}N^{\frac{7\gamma }{8}}
         \\
         & \\
    & \qquad \qquad
         + \Big( \, \frac{h}{d} \, \Big)^{-\frac{1}{4}} Q^{-\frac{1}{4}} N^{\frac{7\gamma }{8}} 
         + \Big( \, \frac{h}{d} \, \Big)^{\frac{1}{12}} Q^{\frac{1}{12}} N^{\frac{11\gamma }{12}}
         + \Big( \, \frac{h}{d} \, \Big)^{-\frac{1}{12}} {Q^{-\frac{1}{12}}} N^{\frac{23\gamma }{24}}  
         \Big) \\
       & \\  
    & \ll
         N^{\varepsilon }
          \left ( D N^{\gamma } Q^{-\frac{1}{2}}
         + D  Q^{\frac{1}{4}} N^{\frac{1}{4} + \frac{5\gamma }{8}}
         + D^{\frac{5}{4}} Q^{-\frac{1}{4}} N^{\frac{7\gamma }{8}} 
         + D Q^{\frac{1}{12}} N^{\frac{1}{12} + \frac{5\gamma }{6}} 
         + D^{\frac{13}{12}} Q^{-\frac{1}{12}} N^{\frac{23\gamma }{24}}  \right ).
\end{align*}
With the choice of $Q$ which we made in \eqref{02.04.2016.eq10} it is now clear that
\begin{equation}\label{104}
 \Omega _{4, 3} \ll \frac{N^{2\gamma - 1}}{(\log N)^{2}}.
\end{equation}

\bigskip

Now, let us carry on with the study of $\Omega_{4, 1}$.
We have choosen the constant $\alpha _{1}$ in such a way, that from \eqref{10.03.2016.eq20} and 
\eqref{15.03.2016.eq10} it follows that
$Y''(l) \asymp v q m N^{-\gamma }$ uniformly for $l \in \left ( L_1, L_2 \right ]$.
Then the sum $\sum_{\frac{P}{m} < l \le \frac{2P}{m} } e \left( Y(l) \right)$
can be bounded 
by the expression on the right side of \eqref{16.03.2016.eq180}. 
This observation illustrates that $\Omega_{4, 1}$ is bounded by the same quantity as $\Omega _{4, 3}$, i.e.
\begin{equation}\label{200}
\Omega _{4, 1} \ll \frac{N^{2\gamma - 1}}{(\log N)^{2}}.
\end{equation}

\bigskip

In a very similar manner one can show that
\begin{equation}\label{201}
\Omega _{4, 2} \ll \frac{N^{2\gamma - 1}}{(\log N)^{2}}.
\end{equation}

\bigskip

Then, from \eqref{132}, \eqref{103} and \eqref{104}  -- \eqref{201} we get
\begin{equation}\label{135}
\Omega_{4} \ll \frac{N^{2\gamma - 1}}{(\log N)^{2}}.
\end{equation}

\bigskip

It remains to find a bound for $\Omega_3$. The same argument as the one for $\Omega_4$ 
can be applied here once more to show that
\begin{equation}\label{202}
\Omega_{3} \ll \frac{N^{2\gamma - 1}}{(\log N)^{2}}.
\end{equation}

\bigskip

From \eqref{136}, \eqref{18.02.eq50}, \eqref{135} and \eqref{202} we conclude that \eqref{90} is satisfied
and the theorem is proved.

\bibliographystyle{plain}

\vspace{5mm}

\noindent Faculty of Mathematics and Informatics \\
Sofia University “St. Kl. Ohridski” \\
5 J.Bourchier, 1164 Sofia, Bulgaria \\
\vspace{5mm}

\noindent {zhpetrov@fmi.uni-sofia.bg \\
dtolev@fmi.uni-sofia.bg}

\end{document}